\documentclass[runningheads]{lncse}
\usepackage{hyperref}
\usepackage{pdfsync}
\usepackage{amsmath}
\usepackage{amsfonts}
\usepackage{amssymb}
\usepackage{epsfig}
\usepackage{graphics} 
\usepackage{xspace}
\usepackage{mathrsfs} 
\usepackage{bbold}    
\usepackage{stmaryrd} 
\usepackage{eufrak}     
\usepackage{subfigure}
\newcommand{\mathscript}{\mathscr}
	   
\providecommand{\figwidth}{\textwidth}
\newcommand{\figscale}{1.1}

\newcommand{\cE}{\ensuremath{\mathscript E}\xspace}

\newcommand{\cL}{\ensuremath{\mathscript L}\xspace}

\newcommand{\cT}{\ensuremath{\mathscript T}\xspace}

\newcommand{\hJ}{\ensuremath{\mathcal J}\xspace}

\newcommand{\hR}{\ensuremath{\mathcal R}\xspace}

\newcommand{\rN}{\ensuremath{\mathbb N}\xspace}

\newcommand{\rP}{\ensuremath{\mathbb P}\xspace}

\newcommand{\rR}{\ensuremath{\mathbb R}\xspace}

\newcommand{\rV}{\ensuremath{\mathbb V}\xspace}
\newcommand{\rW}{\ensuremath{\mathbb W}\xspace}

\newcommand{\Oh} {\operatorname{O}}                   

\newcommand{\naturals}{\rN\xspace}

\newcommand{\reals}{\rR}
\newcommand{\R}[1]{\reals^{#1}}

\newcommand{\lap}{\Delta}
\newcommand{\inflap}{\lap_\infty}
\newcommand{\plap}{\lap_p}

\newcommand{\avg}[1]{\lbrace #1 \rbrace}

\newcommand{\doublelbrace}{
  \mathrel{\ooalign{\cr\kern+1.0pt$\lbrace$\cr\kern-0.5pt$\lbrace$}}}
\newcommand{\doublerbrace}{
  \mathrel{\ooalign{$\rbrace$\cr\kern-1.5pt$\rbrace$\cr\kern-1.0pt}}}

\newcommand{\jump}[1]{\ensuremath{\left\llbracket #1\right\rrbracket}}
\renewcommand{\avg}[1]{\doublelbrace #1 \doublerbrace}
\newcommand{\tjump}[1]{\jump{#1\otimes}}

\providecommand{\qb}[1]{\ensuremath{\!\left[{#1}\right]}}
\newcommand{\qp}[1]{\ensuremath{\!\left({#1}\right)}}
\newcommand{\frob}[2]{\ensuremath{{#1}{:}{#2}}}
\providecommand{\D}{\ensuremath{\mathrm{D}}}
\newcommand{\Hess}{\ensuremath{\D^2}}
\newcommand{\norm}[1]{\ensuremath{\left|#1\right|}}
\newcommand{\Norm}[1]{\ensuremath{\left\|#1\right\|}}
\newcommand{\transpose}{{\boldsymbol\intercal}}   
\newcommand{\Transpose}[1]{\ensuremath{{#1}^{\transpose}}}
\newcommand{\trace}{\operatorname{trace}}             

\newcommand{\tr}[1]{\trace{#1}}
\newcommand{\W}{\ensuremath{\Omega}\xspace}

\newcommand{\CC}{\ensuremath{\operatorname C}\xspace}
\newcommand{\HH}{\ensuremath{\operatorname H}\xspace}
\newcommand{\LL}{\ensuremath{\operatorname L}\xspace}

\newcommand{\WW}{\ensuremath{\operatorname W}\xspace}

\newcommand{\cont}[1]{\ensuremath{\CC^{#1}}}

\newcommand{\leb}[1]{\ensuremath{\LL_{#1}}}

\newcommand{\sob}[2]{\ensuremath{\WW^{#1}_{#2}}}

\newcommand{\sobh}[1]{\ensuremath{\HH^{#1}}}

\newcommand{\Lip}{\ensuremath{\operatorname{Lip}}}

\newcommand{\poly}[1]{\ensuremath{\rP}^{#1}}

\newcommand{\fespace}{\rV}

\newcommand{\fes}{\ensuremath{\fespace}}

\newcommand{\esssup}{\operatorname{ess\,sup}}         

\newcommand{\Forall}{\:\forall\:}

\newcommand{\Foreach}{\quad\Forall}

\newcommand{\closure}[1]{\overline{#1}}

\newcommand{\ie}{i.e.,\xspace}

\renewcommand{\vec}[1]{\ensuremath{\boldsymbol{#1}}}
\newcommand{\geovec}[1]{\ensuremath{\boldsymbol{#1}}}
\newcommand{\geomat}[1]{\ensuremath{\boldsymbol{#1}}}

\newcommand{\pd}[2]{\ensuremath{\partial_{#1}{#2}}\xspace} 

\newcommand{\pdt}[1]{\pd t{#1}}                       
\newcommand{\T}[1]{\cT^{#1}}
\newcommand{\intersected}{\ensuremath{\cap}}
\newcommand{\meet}{\intersected}
\newcommand{\union}[1]{\ensuremath{\bigcup}_{#1}}
\newcommand{\E}{\cE}
\newcommand{\funk}[3]{\ensuremath{#1:#2\to#3}}

\newcommand{\cg}{\ensuremath{\fes_C}}

\newcommand{\dg}{\ensuremath{\fes_D}}
\newcommand{\duality}[2]{\ensuremath{\left\langle #1\,\vert\,#2\right\rangle}}

\renewcommand{\H}{\ensuremath{\geovec{H}}}

\newcommand{\fesW}{\ensuremath{\rW}}
\newcommand{\La}{\ensuremath{\varLambda}\xspace}

\newcommand{\laginterpol}[1]{\La^{#1}}

\newcommand{\A}{\ensuremath{\geovec{A}}}

\renewcommand{\star}[1]{\tilde{#1}}

\newcommand{\residual}{\hR}
\newcommand{\jresidual}{\hJ}

\newcommand{\Program}[1]{\textsf{#1}\xspace}

\newcommand{\fenics}{\Program{FEniCS}\xspace}

\newcommand{\paraview}{\Program{ParaView}\xspace}

\renewcommand{\div}[1]{\operatorname{div}\qp{#1}}

\begin{document}

\title{An adaptive finite element method for the infinity {Laplacian}}

\titlerunning{An AFEM for the Infinity Laplacian}

\author{ Omar Lakkis \inst{1} \and Tristan Pryer \inst{2}   }

\authorrunning{Omar Lakkis and Tristan Pryer}   

\institute{
University of Sussex,
Department of Mathematics,
University of Sussex, 
Brighton,
GB-BN1 9QH
England UK,
\href{mailto:lakkis.o.maths@gmail.com}{\texttt{lakkis.o.maths@gmail.com}}.
\and 
University of Reading,
Department of Mathematics and Statistics,
Whiteknights,
PO Box 220,
Reading,
GB-RG6 6AX,
England UK,
\href{mailto:T.Pryer@Reading.ac.uk}{\texttt{T.Pryer@Reading.ac.uk}}.
}

\maketitle

\abstract{
  We construct a finite element method (FEM) for the infinity
  Laplacian. Solutions of this problem are well
  known to be singular in nature so we have taken the opportunity to
  conduct an a posteriori analysis of the method deriving residual
  based estimators to drive an adaptive algorithm. It is numerically
  shown that optimal convergence rates are regained using the adaptive
  procedure.
}

\section{Introduction}

\label{sec:intro}

Nonlinear partial differential equations (PDEs) arise in many
areas. Their numerical simulation is extremely important due to the
additional difficulties arising in their classical solution
\cite{CaffarelliCabre:1995}. One such example is that of the \emph{infinity
Laplace operator} $\inflap$ defined by
\begin{equation}
  \label{eq:inf-lap}
  \inflap u 
  :=
  \frac
    {\sum_{i=1}^d\sum_{j=1}^d 
  \partial_i u \partial_j u
  \partial_{{i}{j}} u}
    {\sum_{i=1}^d\qp{\partial_i u}^2}
  =
  \frac{\frob{\qp{\nabla u \otimes \nabla u}}{\Hess u}}{\norm{\nabla u}^2},
\end{equation}
for a twice-differentiable function $\funk u\W\reals$, $\W\in\R d$
open, bounded and connected, where
\begin{equation}
  \nabla u := 
  \begin{bmatrix}
    \partial_1 u\\\vdots\\\partial_d u
  \end{bmatrix},
  \quad
  \vec x\otimes\vec y:= \vec x\Transpose{\vec y},
  \text{ and }
  \frob{\vec X}{\vec Y}:=\tr{\Transpose{\vec X} \vec Y}
\end{equation}
denote, respectively, the gradient, the (algebraic) tensor product of
$\vec x,\vec y\in\R d$, and the Frobenius inner product of two
matrices $\vec X,\vec Y\in\R{d\times d}$.  This equation has been
popular in classical studies
\cite[e.g.]{Aronsson:1986,BarronEvansJensen:2008} but is difficult to
pose numerical schemes due to its nondivergence structure and general
lack of classical solvability.  The infinity Laplacian, which is in
fact a misnomer (\emph{homogeneous infinity Laplacian} is more
precise), occurs as the weighted formal limit of a variational
problem. A more appropriate terminology would be that of
\emph{infinite harmonic} function $u$ being one that solves $\inflap
u=0$.  This is justified, at least heuristically, as being the formal
limit of the $p$-harmonic functions, $u_p$, $p\geq1$, $p\to\infty$
where
\begin{equation}
  \label{eq:p-lap-to-inf-lap}
  \begin{split}
    0
    =
    \plap u_p
    :=
    \div{\norm{\nabla u_p}^{p-2}\nabla u_p}
    = 
    \norm{\nabla u_p}^{p-2} \lap u_p
    +
    \qp{p-2}\norm{\nabla u_p}^{p-2} \inflap u_p.
  \end{split}
\end{equation}
Multiplying by $\norm{\nabla u_p}^{2-p}/(p-2)$ and taking the limit as
$p\to\infty$ it follows that a would be limit $u=\lim_{p\to\infty}u_p$
is infinite harmonic. A rigorous treatment is provided in
\cite{CrandallEvansGariepy:2001} and is based on the variational
observation that the Dirichlet problem for the $p$--Laplacian is the
Euler--Lagrange equation of the following \emph{energy} functional
\begin{equation}
  \cL_p[u]
  := 
  \frac{1}{p}\Norm{u}_{\leb{p}(\W)}^p 
  =
  \int_\W \frac{1}{p} \norm{\nabla u}^p \text{ for } p \in [1,\infty)
\end{equation}
with appropriate Dirichlet boundary conditions.  By analogy, setting
\begin{equation}
  \label{eq:infty-Dirichlet-integral}
  \cL_\infty[u]
  :=
  \Norm{\nabla u}_{\leb{\infty}(\W)}
  = 
  \esssup_{\W} \norm{\nabla u},
\end{equation}
we seek $u\in\Lip(\W)=\sob1\infty(\W)$, the space of Lipschitz continuous functions
over $\W$ (Rademacher), with $u = g$ on $\partial \W$ such that
\begin{equation}
  \cL_\infty[u] \leq \cL_\infty[v] \Foreach v\in\Lip(\W) 
  \text{ and } v=g \text{ on }  \partial\W.
\end{equation}
Show that the solution exists and define it to be infinite harmonic.
Such a solution is called \emph{absolutely minimising Lipschitz
  extension of $g$}, we call it infinite harmonic.  The infinity
Laplacian is thus considered to be the paradigm of a variational
problem in $\sob1{\infty}(\W)$.

If the solution is smooth, say in $\cont2$ and has no internal
extrema, it can be shown to satisfy (\ref{eq:p-lap-to-inf-lap})
classically.  But an infinite harmonic function is generally not a
classical solution (those in $\cont{2}(\W)$ satisfying
(\ref{eq:inf-lap}) everywhere.  Therefore solutions of
(\ref{eq:p-lap-to-inf-lap}) must be sought in a weaker sense. The
notion of viscosity solution, introduced for second order PDEs in
\cite{CrandallIshiiLions:1992} turns out to be the correct setting to
seek weaker solutions.  Existence and uniqueness of a viscosity
solution to the homogeneous infinity Laplacian (\ref{eq:inf-lap}) has
been studied \cite{Jensen:1993}.  If the domain $\W$ is bounded, open
and connected then (\ref{eq:inf-lap}) has a unique viscosity solution
$u \in \cont{0}(\closure{\W})$.  In the case $\W\subset\reals^2$ this
can be improved to $u\in\cont{1,\alpha}(\closure\W)$
\cite{EvansSavin:2008}.  A study of existence and uniqueness of
viscosity solutions to the inhomogeneous infinity Laplacian can be
found in \cite{LuWang:2008}. With $\W$ defined as before and in
addition if $f\in\cont{0}(\W)$ and does not change sign, \ie $\inf_\W
f > 0$ or $\sup_\W f < 0$, one can find a unique viscosity solution.

As to the topic of numerical methods to approximate the infinity
Laplacian, to the authors knowledge only two methods exist. The first
is based on finite differences \cite{Oberman:2005}. The scheme
involves constructing monotone sequences of schemes over concurrent
lattices by minimising the discrete Lipschitz constant over each node
of the lattice. The second is a finite element scheme named the
vanishing moment method \cite{FengNeilan:2009} in which the 2nd order
nonlinear PDE is approximated via sequences of biharmonic quasilinear
4th order PDEs.

In this paper we present a finite element method for the infinity
Laplacian, without having to deal with the added complications of
approximating a 4th order operator. It is based on the non-variational
finite element method introduced in \cite{LakkisPryer:2011a}. Roughly,
this method involves representing the \emph{finite element Hessian}
(see Definition \ref{defn:feh}) as an auxiliary variable in the
formulation, to deal with the nonvariational structure. We also
consider the problem as the steady state of an evolution equation
making use of a \emph{Laplacian relaxation} technique (see Remark
\ref{rem:evolution}) \cite{Awanou:2011,EsedogluOberman:2011} to
circumvent the degeneracy of the problem.

The structure of the paper is as follows: In
\S\ref{sec:discretisation} we examine the linearisation of the PDE and
present the necessary framework for the discretisation and state an a
posteriori error indicator for the discrete problem. The estimator is
of residual type and is used to drive an adaptive algorithm which is
studied and used for numerical experimentation is
\S\ref{sec:numerics}. We choose our simulations in such a way that
they can be compared with those given in \cite{Oberman:2005,FengNeilan:2009}.

\section{Notation, linearisation and discretisation}
\label{sec:discretisation}

We consider the inhomogeneous Infinity Laplace problem with Dirichlet
boundary conditions on a domain $\W\subset \reals^d$.
\begin{gather}
  \label{eq:model-problem-with-bcs}
  \inflap u = f \quad \text{ in } \W \quad \text{ and }
  \quad 
  u = g \quad \text{ on } \partial \W
\end{gather}
with problem data $f,g\in\cont{0}(\W)$ chosen such that $f$ does not
change sign throughout $\W$. In this case there exists a unique
viscosity solution to (\ref{eq:model-problem-with-bcs})
\cite{LuWang:2008}.
 
\subsubsection{Linearisation of the continuous problem (\ref{eq:inf-lap})}

The application of a standard fixed point linearisation to
(\ref{eq:model-problem-with-bcs}) results in the following sequence of
linear non-divergence PDEs: Given an initial guess $u^0$, for each
$n\in\naturals$ find $u^{n+1}$ such that
\begin{equation}
  \label{eq:continuous-linearised}
  \frob{\frac{\qp{\nabla u^n \otimes \nabla u^n}}{\norm{\nabla u^n}^2}}{\Hess u^{n+1}} = f.
\end{equation}

Due to the degeneracy of the problem we introduce a slightly modified
problem which utilises \emph{Laplacian relaxation} \cite{Awanou:2011,EsedogluOberman:2011},
the problem is to find $u^{n+1}$ such that
\begin{equation}
  \label{eq:continuous-laplace-precon}
  \frob{\qp{\frac{\nabla u^n\otimes \nabla u^n}{\norm{\nabla u^n}^2} + \frac{\geovec I}{\tau}}}{\Hess u^{n+1}} = f + \frac{\Delta u^n}{\tau}
\end{equation}
where $\tau\in\reals^+$.

\begin{remark}
  \label{rem:evolution}
  The discretisation proposed in (\ref{eq:continuous-laplace-precon}) is nothing but an implicit one stage discretisation of the following evolution equation
  \begin{equation}
    \pdt{\qp{\Delta u}} + \inflap u = f,
  \end{equation}
  where $\Delta u$ is used as shorthand for $\Delta_2 u$, the 2--Laplacian.

  With that in mind we must take care with our choice of $\tau$ which
  can be regarded as a timestep. We require a $\tau$ that is large
  enough to guarantee reaching the steady state and small enough such
  that we do not encounter stability problems.
\end{remark}


\subsubsection{Discretisation of the sequence of linear PDEs (\ref{eq:continuous-laplace-precon})}

Let $\T{}$ be a conforming, shape regular triangulation of $\W$, namely, $\T{}$ is a
finite family of sets such that
\begin{enumerate}
\item $K\in\T{}$ implies $K$ is an open simplex (segment for $d=1$,
  triangle for $d=2$, tetrahedron for $d=3$),
\item for any $K,J\in\T{}$ we have that $\closure K\meet\closure J$ is
  a full sub-simplex (i.e., it is either $\emptyset$, a vertex, an
  edge, a face, or the whole of $\closure K$ and $\closure J$) of both
  $\closure K$ and $\closure J$ and
\item $\union{K\in\T{}}\closure K=\closure\W$.
\end{enumerate}
We also define $\E$ to be the skeleton of the triangulation, that is
the set of sub-simplexes of $\T{}$ contained in $\W$ but not $\partial
\W$. For $d=2$, for example, $\E$ would consist of the set of edges of
$\T{}$ not on the boundary. We also use the convention where
$h(\vec{x}):=\max_{\closure K\ni \vec x}h_K$ to be the mesh-size
function of $\T{}$.

 \begin{definition}[continuous and discontinuous FE spaces]
   Let $\poly k(\T{})$ denote the space of piecewise polynomials of
   degree $k$ over the triangulation $\T{}$ of $\W$. We introduce the \emph{finite element spaces}
   \begin{gather}
     \label{eqn:def:finite-element-space}
     \dg\qp{k} = \poly k(\T{})
     \qquad 
     \cg\qp{k} = \poly k(\T{}) \cap \cont{0}(\W)
   \end{gather}
   to be the usual spaces of discontinuous and continuous piecewise polynomial
   functions over $\W$.
 \end{definition}
 
\begin{remark}[generalised Hessian]
  \label{rem:generalised-hessian}
  Given a function $v\in\sobh1(\W)$ and let $\geovec
  n:\partial\W\to\reals^d$ be the outward pointing normal of $\W$ then the
  \emph{generalised Hessian} of $v$, $\Hess v$ satisfies the following identity:
  \begin{equation}
    \label{eq:generalised-hessian}
    \duality{\Hess v}{\phi}
    = 
    -
    \int_\W {\nabla v}\otimes{\nabla \phi}
    +
    \int_{\partial\W}{\nabla v} \otimes {\geovec n \ \phi} 
    \Foreach \phi\in\sobh1(\W),
  \end{equation}
  where the final term is understood as a duality pairing between
  $\sobh{-1/2}(\partial\W) \times \sobh{1/2}(\partial\W)$.
\end{remark}

\begin{remark}[nonconforming generalised Hessian]
  \label{rem:nonconforming-generalised-hessian}
  The test functions applied to define the generalised Hessian in
  Remark \ref{rem:generalised-hessian} need not be
  $\sobh{1}(\W)$. Suppose they are $\sobh{1}(K)$ for each $K\in\T{}$
  then it is clear that
  \begin{equation}
    \begin{split}
      \duality{\Hess v}{\phi} 
      &=
      \sum_{K\in\T{}} 
      \qp{-\int_K{\nabla v}\otimes{\nabla \phi} 
      +
      \int_{\partial K}{\nabla v}\otimes{\geovec n_K \phi}}
      \\
      &=
      \sum_{K\in\T{}} 
      -\int_K{\nabla v}\otimes {\nabla \phi} 
      +
      \sum_{e\in\E} 
      \int_e{\avg{\nabla v}}\otimes{\jump{\phi}}
      +
      \sum_{e\in\partial\W} 
      \int_e{\nabla v}\otimes{\geovec n \ \phi},
    \end{split}
  \end{equation}
  where $\jump{\cdot}$ and $\avg{\cdot}$ denote the \emph{jump} and
  \emph{average}, respectively, over an element edge, that is, suppose $e$ is a
  $\qp{d-1}$ subsimplex shared by two elements $K^+$ and $K^-$ with
  outward pointing normals $\geovec n^+$ and $\geovec n^-$
  respectively, then 
  \begin{equation}
    \jump{\geovec \eta}
    =
    {\eta\big|_{K^+}} \geovec n^+
    +
    {\eta\big|_{K^-}} \geovec n^-
    \text{ and }
    \avg{\geovec \xi}
    = 
    \frac{1}{2}
    \qp{\geovec \xi\big|_{K^+} + \geovec\xi\big|_{K^-}}.
  \end{equation}
\end{remark}

\begin{definition}[finite element Hessian]
  \label{defn:feh}
  From Remark \ref{rem:generalised-hessian} and Remark
  \ref{rem:nonconforming-generalised-hessian} for $V\in\cg\qp{k}$ we
  define the \emph{finite element Hessian}, $\H[V]\in\qb{\dg\qp{k}}^{d\times
    d}$ such that we have
  \begin{equation}
    \int_\W{\H[V]}{\phi} 
    =
    \duality{\Hess V}{\phi} \Foreach \phi\in\dg\qp{k}.
  \end{equation}
\end{definition}

We discretise (\ref{eq:continuous-laplace-precon}) utilising the
non-variational Galerkin procedure proposed in
\cite{LakkisPryer:2011a}. We construct finite element spaces
$\fes:=\cg\qp{k}$ and $\fesW$ which can be taken as $\cg\qp{k}$,
$\dg\qp{k}$ or $\dg\qp{k-1}$. Then given ${U^0} = \laginterpol{} u^0$,
for each $n\in\naturals_0$ we seek $\qp{U^{n+1}, \H[U^{n+1}]} \in \fes
\times \qb{\fesW}^{d\times d}$ such that
\begin{equation}
  \begin{split}
  \label{eq:2-0-NVFEM}
  &
  \int_\W{\frob{\qp{\frac{\nabla U^n \otimes \nabla U^n}{\norm{\nabla
            U^n}^2} + \frac{\geovec I}{\tau}}}{\H[U^{n+1}]}}{\Psi} =
  \int_\W\qp{f + \frac{\trace{\H[U^n]}}{\tau}}{\Psi}
  \\
  &\int_\W{\H[U^{n+1}]}{\Phi} = - \int_\W{\nabla U^{n+1}}\otimes
  {\nabla \Phi} + \sum_{e\in\E} \int_e{\avg{\nabla
      U^{n+1}}}\otimes{\jump{\Phi}}
  \\
  &\qquad\qquad\qquad\qquad\qquad\qquad+ \sum_{e\in\partial\W}
  \int_e{\nabla U^{n+1}}\otimes{\geovec n \ \Phi} \Foreach \qp{\Psi,
    \Phi} \in \fes\times\fesW.
  \end{split}
\end{equation}

\begin{remark}[computational efficiency]
  Making use of a $\dg\qp{k}$ or $\dg\qp{k-1}$ space to represent the
  finite element Hessian allows us to construct a much faster
  algorithm in comparison to using a $\cg{k}$ space for $\fesW$ due to
  the local representation of the $\leb{2}(\W)$ projection of
  discontinuous spaces \cite{DednerPryer:2012}.
\end{remark}

\begin{theorem}[a posteriori residual upper error bound]
  \label{the:aposteriori}
  Let $u$ be the solution to the infinity Laplacian
  \eqref{eq:model-problem-with-bcs} and $U^n$ be the $n$-th step in
  the linearisation defined by \eqref{eq:2-0-NVFEM}. Let
  \begin{equation}
    \A[v] := \frac{\nabla v \otimes \nabla v}{\norm{\nabla v}^2} + \frac{\geomat I}{\tau},
  \end{equation}
  then there exists a $C > 0$ such that
  \begin{equation}
    \begin{split}
    \label{eq:aposteriori-estimate}
    \Norm{f + \frac{\Delta U^n}{\tau} - \frob{\A[U^n]}{\Hess U^{n+1}}}_{\sobh{-1}(\W)} 
    &\leq
    C\bigg( 
    \sum_{K\in\T{}} 
    h_K \Norm{\residual[U^{n},U^{n+1},f]}_{\leb{2}(K)}
    \\
    & \qquad\qquad+    
    \sum_{e\in\E}
    h_K^{1/2} \Norm{\jresidual[U^{n},U^{n+1}]}_{\leb{2}(e)}
    \bigg)
    \end{split}
  \end{equation}
  where the interior residual, $\residual[U,\A,f]$, over a simplex
  $K$ and jump residual, $\jresidual[U,\A]$, over a common wall $e =
  \closure{K}^+\cap\closure{K}^-$ of two simplexes, $K^+$ and $K^-$
  are defined as
  \begin{gather}
    \Norm{\residual[U^n, U^{n+1},f]}^2_{\leb{2}(K)} 
    =
    \int_K\qp{{f - \frob{\A[U^n]}{\Hess U^{n+1}}} + \frac{\Delta U^n}{\tau}}^2,
    \\
    \Norm{\jresidual[U^n,U^{n+1}]}^2_{\leb{2}(e)}
    =
    \int_e \qp{
      \frac{\jump{\nabla U^n}}{\tau} - \frob{\A[U^n]}{\tjump{\nabla U^{n+1}}}
    }^2,
  \end{gather}
  with 
  \begin{equation}
    \tjump{\geovec \xi} := \geovec\xi|_{K^+} \otimes \geovec n^+ + \geovec\xi|_{K^-} \otimes \geovec n^-,
  \end{equation}
  being defined as a \emph{tensor jump}.
\end{theorem}

\section{Numerical experiments}
\label{sec:numerics}

All of the numerical experiments in this section are implemented using
\fenics and visualised with \paraview. Each of the tests are on the
domain $\W = [-1,1]^2$, choosing the finite element spaces $\fes =
\cg\qp{1}$ and $\fesW = \dg\qp{0}$. This is computationally the
quickest implementation of the non-variational finite element method
and the lowest order stable pair of FE spaces for this class of
problem.


\subsubsection{Benchmarking and convergence -- Classical solution}

\begin{figure}[t]
  \caption[]
          {
            \label{fig:convergence-classical}
            We benchmark the approximation of a classical solution to
            the inhomogeneous infinity Laplacian, plotting the log of
            the error together with its estimated order of
            convergence. We examine both $\leb{2}(\W)$ and
            $\sobh{1}(\W)$ norms of the error together with the
            residual estimator given in Theorem
            \ref{the:aposteriori}. The linearisation tolerance is
            coupled to the mesh-size such that the linearisation is run
            until $\Norm{U^n - U^{n-1}}\leq 10 h^2$. The convergence
            rates are optimal, that is, $\Norm{u - U^N} = \Oh(h^2)$
            and $\norm{u - U^N}_1 = \Oh(h)$.  }
  \begin{center}
    \subfigure[][{Convergence rates}]{
      \includegraphics[scale=\figscale, width=0.47\figwidth]{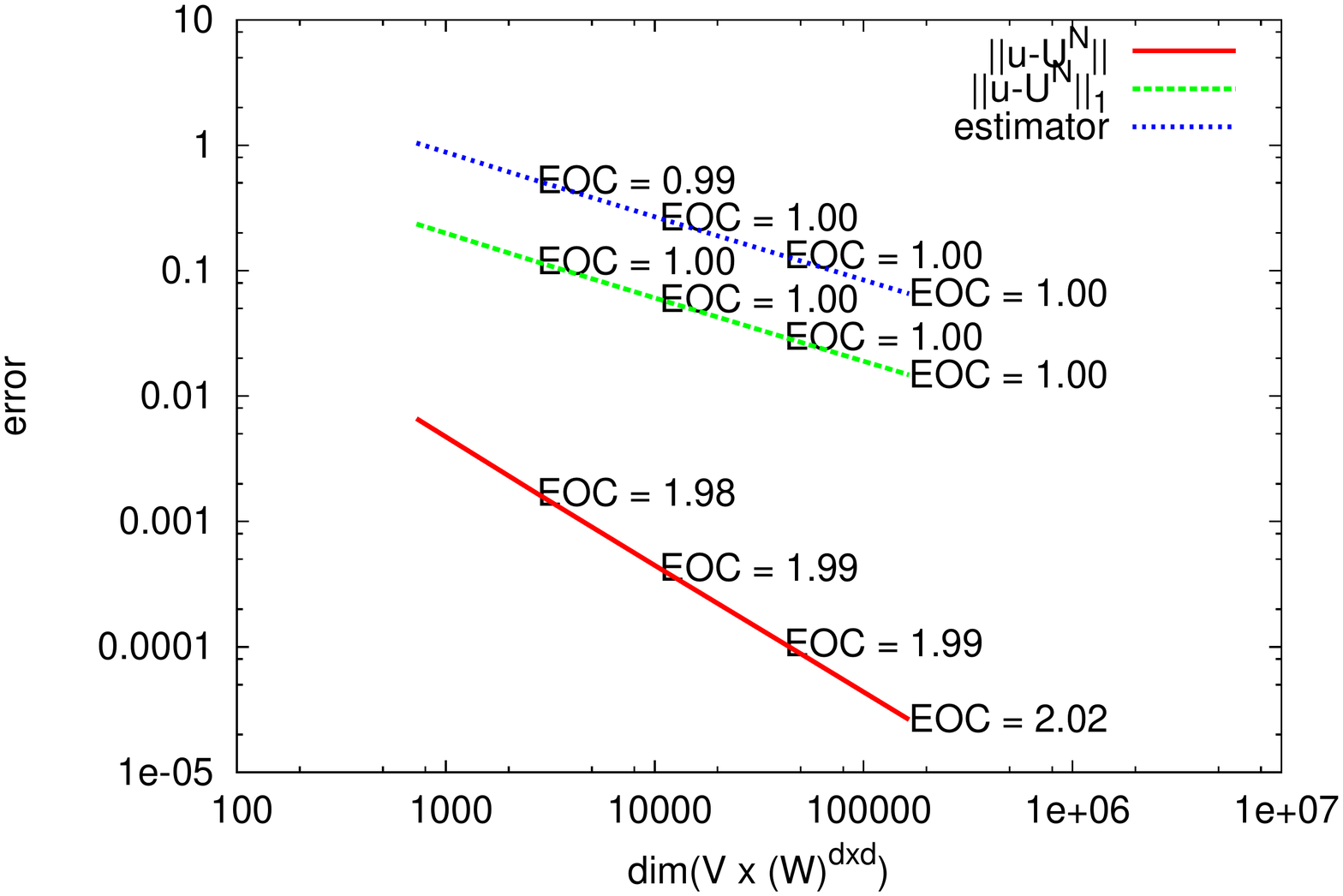}
    }
    \subfigure[][{Finite element approximation}]{
      \includegraphics[scale=\figscale, width=0.47\figwidth]{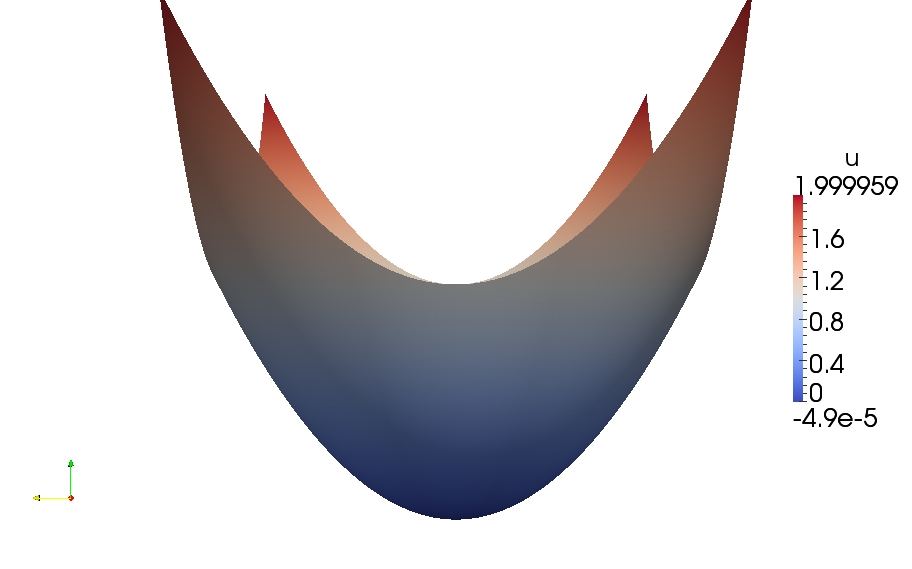}
    }
  \end{center}
\end{figure}

To benchmark the numerical algorithm we choose the data $f$ and $g$
such that the solution is known and classical. In the first instance
we choose
$f \equiv 2$
and
$g = \norm{\geovec x}^2$.
It is easily verified that the exact solution is given by
$u = \norm{\geovec x}^2$.
Figure \ref{fig:convergence-classical} details a numerical experiment
on this problem.

\begin{remark}[on the value of $\tau$]
  The optimal values of the \emph{timestep parameter} or tuning
  parameter $\tau$ depend upon the regularity of the solution. For
  example, for a classical solution, one may
  choose $\tau$ large. In the numerical experiment above we took
  $\tau = 1000$. Since the linearisation is nothing more than seeking
  the steady state of the evolution equation
  \eqref{eq:continuous-laplace-precon}. The convergence (in $n$) is
  extremely quick taking no more than five iterations. 

  For the examples below one must be careful choosing $\tau$, we will
  be looking at viscosity solutions that are not $\cont{2}(\W)$, in
  this case the lack of regularity of the solution will lead to an
  unstable linearisation for large $\tau$. In each of the cases below
  $\tau \in [1:10]$ was sufficiently small to achieve convergence of
  the linearisation in at most twenty iterations. 
\end{remark}

\subsubsection{A known viscosity solution to the homogeneous problem}

\begin{figure}[t]
  \caption[]
          {
            \label{fig:convergence-arronson}
            We benchmark problem (\ref{prob:viscocity-homog}),
            plotting the log of the error together with its estimated
            order of convergence. We examine both $\leb{2}(\W)$ and
            $\sobh{1}(\W)$ norms of the error together with the
            residual estimator given in Theorem
            \ref{the:aposteriori}. We choose $\tau = 1$ and
            the linearisation tolerance is coupled to the mesh-size such
            that the linearisation is run until $\Norm{U^n -
              U^{n-1}}\leq 10 h^2$. The convergence rates are
            suboptimal due to the singularity, that is, $\Norm{u - U^N}
            \approx \Oh(h^{1.8})$ and $\norm{u - U^N}_1 \approx
            \Oh(h^{0.8})$.  }
  \begin{center}
    \subfigure[][{Convergence rates}]{
      \includegraphics[scale=\figscale, width=0.47\figwidth]{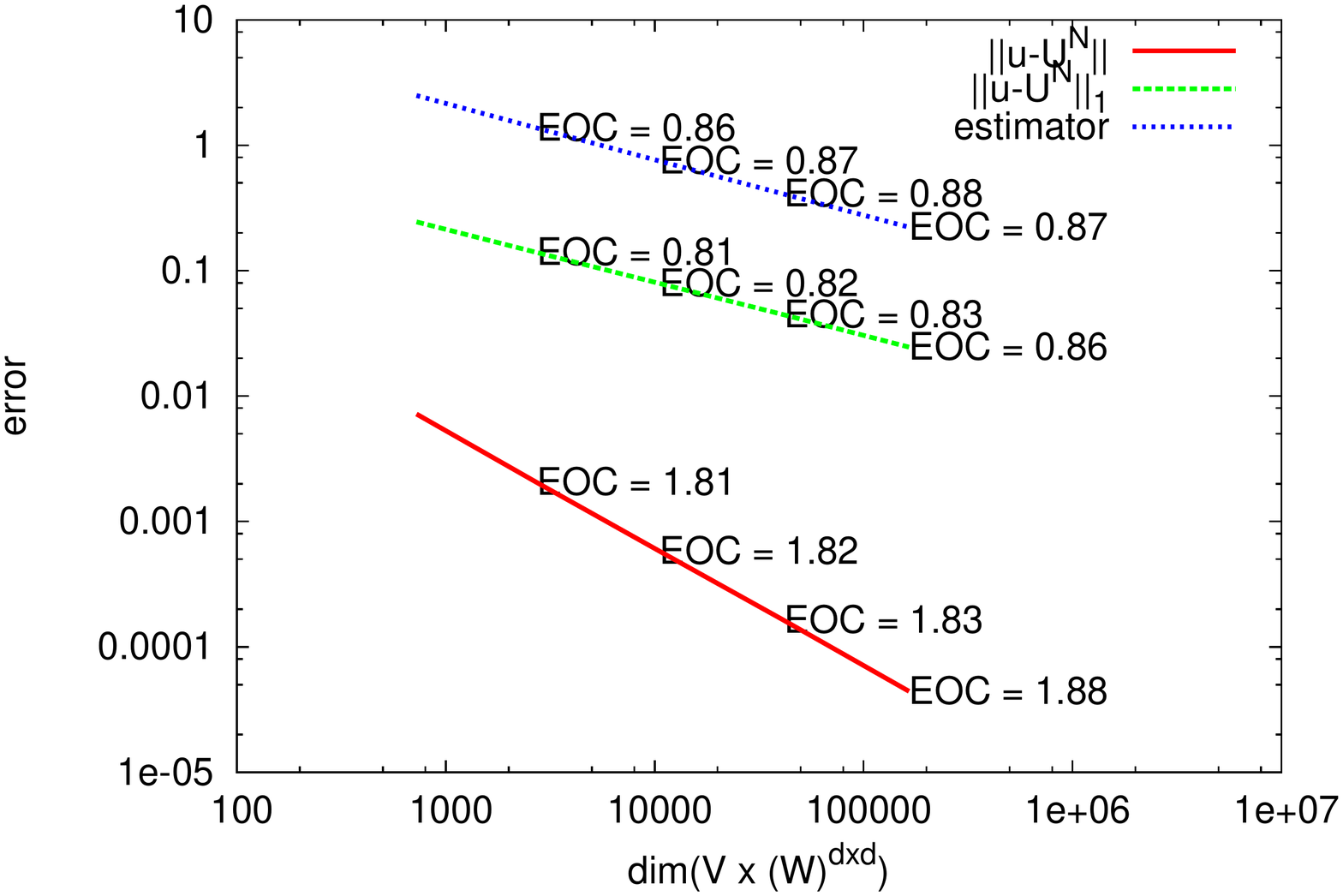}
    }
    \subfigure[][{Finite element approximation}]{
      \includegraphics[scale=\figscale, width=0.47\figwidth]{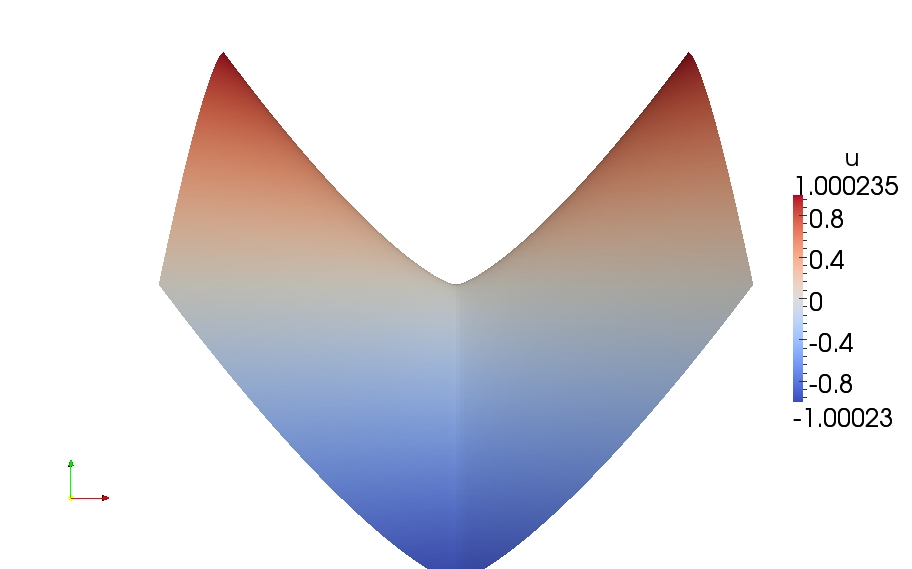}
    }
  \end{center}
\end{figure}

To test the convergence of the method applied to a singular solution
of the homogeneous problem we fix
\begin{gather}
  \label{prob:viscocity-homog}
  f \equiv 0
  \text{ and }
  g = \norm{x}^{4/3} - \norm{y}^{4/3},
\end{gather}
where $\geovec x = \Transpose{\qp{x,y}}$. A known viscosity solution
of this equation is the Aronsson solution \cite{Aronsson:1986},
\begin{equation}
  u(\geovec x) = \norm{x}^{4/3} - \norm{y}^{4/3}.
\end{equation}
The function has singular derivatives about the coordinate axis, in
fact $u\in\cont{1,1/3}(\W)$. Figure \ref{fig:convergence-arronson}
details a numerical experiment on this problem.

\begin{figure}[ht]
  \caption[]
          {
            \label{fig:adapt-arron-solution}
            This is an adaptive approximation of the viscosity
            solution $u = \norm{x}^{4/3} - \norm{y}^{4/3}$ from
            (\ref{prob:viscocity-homog}). The estimator tolerance was set
            at $0.1$ to coincide with the final estimate from the
            benchmark solution from Figure
            \ref{fig:convergence-arronson}. The final number of
            degrees of freedom was $36,325$ compared to the uniform scheme
            which took $165,125$ degrees of freedom to reach the same
            tolerance. We chose $\tau = 0.1$ as the timestep
            parameter.  }
  \begin{center}
    \subfigure[][{The finite element approximation viewed from the top.}]{
      \includegraphics[scale=\figscale, width=0.47\figwidth]{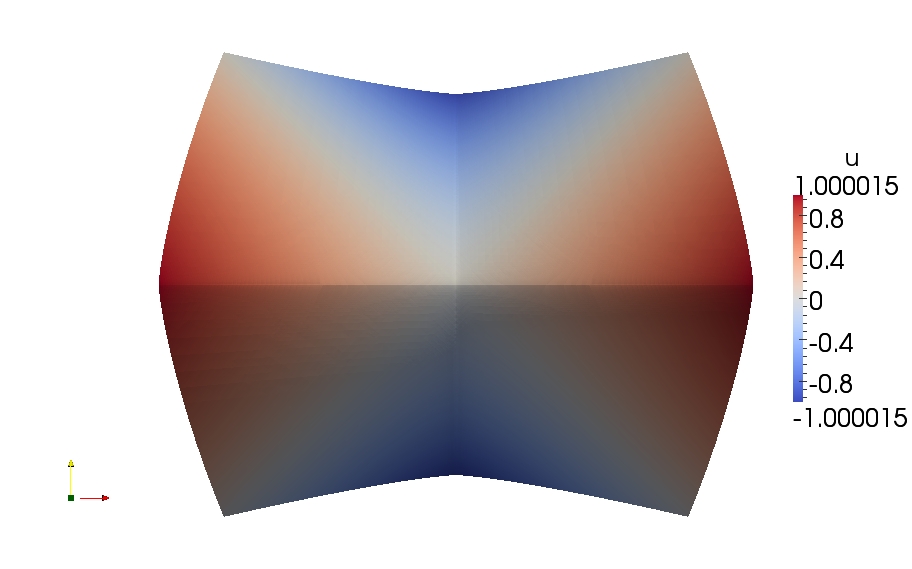}
    }
    \hfill
    \subfigure[][{The underlying mesh}]{
      \includegraphics[scale=\figscale, width=0.47\figwidth]{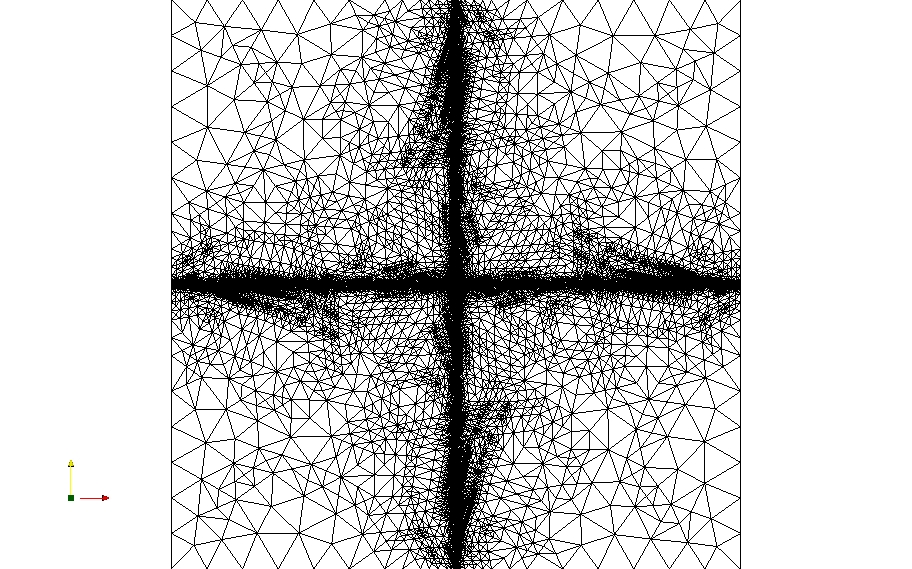}
    }
  \end{center}
\end{figure}

In Figure \ref{fig:adapt-arron-solution} we conduct an adaptive
experiment based on the newest vertex bisection method. 

\ifx\undefined\bysame
\newcommand{\bysame}{\leavevmode\hbox to3em{\hrulefill}\,}
\fi

\end{document}